\documentclass[12pt]{article}
\usepackage[utf8]{inputenc}  \usepackage[T1]{fontenc}  \usepackage{lmodern}
\usepackage{array} 
\usepackage{paralist} 
\usepackage{amsmath} \usepackage{amsthm} 
\usepackage[british]{babel}  \usepackage{amssymb}  \usepackage{graphicx}
\usepackage[left=2.5cm,top=2.5cm,right=2.5cm,bottom=2cm]{geometry}  \usepackage{color} 
\usepackage{comment}

\newtheorem {theorem} {Theorem}

\numberwithin{equation}{section}
\newcommand{\R}{\mathbb{R}}

\title{Limit cycles of linear vector fields on $(\mathbb{S}^2)^m \times \R^ n$}

\date{}

\begin{document}
	\maketitle

	\vspace{-1.5cm}
	
	\centerline{Clara Cufí-Cabré and Jaume Llibre}
	\medskip
	{\footnotesize
		\centerline{Departament de Matemàtiques} 
		\centerline{Universitat Autònoma de Barcelona (UAB)} 
		\centerline{08193 Bellaterra, Barcelona, Spain}
		\centerline{\texttt{clara@mat.uab.cat, \, jllibre@mat.uab.cat}}
	} 

	\bigskip

	\begin{abstract}
		It  is well known  that  linear vector fields defined in $\R^n$ can not have limit  cycles, but this is not the case for linear  vector fields defined in  other manifolds. We study the existence of limit cycles bifurcating from a continuum of periodic orbits of linear vector fields on manifolds of the form $(\mathbb{S}^2)^m \times \R^ n$ when such vector fields are perturbed inside the class of all linear vector fields. The study is done using the averaging theory. We also present  an open problem concerning  the maximum number of limit cycles of linear  vector fields  on  $(\mathbb{S}^2)^m \times \R^ n$.
	\end{abstract}
	
	\vspace{0.1cm}
	
	\begin{footnotesize}
	\end{footnotesize}

	\section{Introduction and statement of the main results}
	
	The study of periodic orbits of differential systems  play an important role in the qualitative theory of ordinary differential equations and their applications. A \emph{limit cycle} is defined as a periodic orbit of a differential system which is isolated in the set of all periodic orbits of the system.  There are many works concerning the study of limit cycles and their applications (see for instance \cite{ChLlo, GiLliVia, HanLi, Ilya} and the references quoted therein).

	It is well know that linear vector fields in $\R^n$ can not have limit cycles, but this is not the case if one considers linear vector fields in other manifolds different from $\R^n$. The objective of this paper is to study the existence of limit cycles of linear vector fields defined on the manifolds $(\mathbb{S}^2)^n \times \R^n$. 
	
	The problem of studying limit cycles of linear vector fields on manifolds different form $\R^n$ was already treated in \cite{lli-zhang}, where the authors consider linear vector fields on  $\mathbb{S}^m \times \R^n$, and they conjecture that such vector fields  may have at most one limit cycle.

	Linear autonomous differential systems, namely, systems of the form $\dot x = A x + b$, where $A$ is a $n \times n$ real matrix and $b$ is a vector in $\R^n$,  are the easiest  systems to study because their solutions can be completely determined (see \cite{arnold, sotomayor}), but still they play an important role in the theory of differential systems. Thus when a nonlinear differential system has a hyperbolic equilibrium point, the dynamics around that point is determined by the linearization of the vector field at that point (Hartman-Grossman theorem, see \cite{hartman}). 
	
	Also linear vector fields having invariant subspaces of periodic orbits can be perturbed inside a concrete class of nonlinear  differential systems to obtain limit cycles of these nonlinear systems bifurcating from  the periodic orbits of the linear  system (see \cite{fer-lli-teix, lli-mart-teix, lli-teix, lli-teix-torr}).  
	
	Moreover,  linear differential systems of the form $\dot x = A x + B u$, where $x$ are the state variables and $u$ is the control input, are applied in control theory for the modeling of hybrid systems (see \cite{laf-1,laf-2}).
	
	The preceeding exemples show some of the importance of linear differential systems. In this paper we show that linear differential systems can have limit cycles when the manifold where they are defined is different from $\R^n$, and we consider the question of how many limit cycles can have at most a linear vector field depending on the manifold where it is defined.
	
	Let $M$ be a smooth connected manifold of dimension $n$, and let $TM$ be its tangent bundle.  A \emph{vector field} on $M$ is a map $X : M \to TM$ such that $X(x) \in T_xM$, where $T_xM$ is the tangent space of $M$ at the point $x$. 
	
	A \emph{linear vector field} in $\R^n$ is a vector field of the form $X(x) = A x + b$, with $x, b \in \R^n$ and where $A$ is a $n \times n$ real matrix.  As it is well known linear vector fields on $\R^n$ either do not have periodic orbits or their periodic orbits form a continuum, and therefore they do not have limit cycles.
	
	In this paper we consider linear vector fields on some manifolds of the form $(\mathbb{S}^2)^m \times \R^n$, where $\mathbb{S}^2$ denotes the unit two-dimensional sphere. Here the sphere $\mathbb{S}^2$ is parameterized by the coordinates $(\theta, \varphi)$, where $\theta \in [-\pi, \pi)$ denotes the azimuth angle and $\varphi \in [-\pi /2,  \pi/2]$ is the polar angle. Hence the curve $\{\varphi = 0 \}$ is the equator of the sphere.
	
	Let $(\theta_1, \varphi_1, \dots \theta_m, \varphi_m,   x_1, \cdots, x_n )$ denote the coordinates of the space $(\mathbb{S}^2)^m \times \R^n$. Then we say that a vector field $X$ is \emph{linear} on $M=(\mathbb{S}^2)^m \times \R^n$ if the expression of $X$ in the coordinates $z = (\theta_1, \varphi_1, \dots \theta_m, \varphi_m,  x_1, \cdots, x_n ) \in M$ is of the form 
	$X(z) = Az + b$, with $b \in M$ and where $A$ is a $(2m + n) \times (2m + n)$ real matrix. 
	
	A simple example in which a linear differential system on the manifold $(\mathbb{S}^2)^m \times \R^n$  has a limit cycle is the following. Take $m=1, \, n=0$ and consider the linear system on the sphere $\mathbb{S}^2$ given by
	$$
	\dot  \theta = 1, \qquad  \dot  \varphi = \varphi,
	$$
	for $\theta \in [-\pi, \pi)$ and $\varphi \in (-\pi/2, \pi/2)$, and
	$$
	\dot \theta = 0, \qquad \dot \varphi = 0
	$$
	for $\varphi = \pm \pi /2$. Then, clearly the equator of the sphere $\{\varphi = 0 \}$ is the only periodic orbit of the  system, and therefore it is a limit cycle. 
	
	In this paper we consider generic linear perturbations of some linear vector fields on three different manifolds of the form $(\mathbb{S}^2)^m \times \R^n$, and we study whether those families of linear differential systems can  have limit cycles.

	Let $M = \mathbb{S}^2 \times \R$ and consider the linear differential system in $M$ given by
	\begin{equation} \label{sist-lineal-1}
		\dot \theta = 1,  \quad \dot \varphi = 0, \quad
		\dot r = r-1,
	\end{equation}
	for $r \in \R$, $\theta \in [-\pi,  \pi)$ and $\varphi \in (-\pi/2, \pi/2)$, and with $\dot \theta = 0$ on  the straight lines $R_1 = \{ \varphi=-\pi/2\}$ and $R_2 = \{ \varphi=\pi/2 \}$.
	
	The solution of system \eqref{sist-lineal-1} is given by
	$$
	\theta(t) = \theta_0 + t, \quad \varphi(t) = \varphi_0, \quad r(t) = (r_0 -1) e^t +1.
	$$
	Thus the sphere $\{r=1\}$ is an invariant manifold with two equilibrium points at the north and the south poles, and is foliated by periodic orbits of period $2 \pi$, corresponding to the parallels of the sphere, except at the poles. Moreover the straight lines $R_1$ and $R_2$ are invariant. 
	
	First we shall study the bifurcation of limit cycles when we perturb system \eqref{sist-lineal-1} inside the class of all linear differential systems, and we shall  see that one of the periodic orbits contained in the sphere $\{r=1\}$ may bifurcate to a limit cycle under certain hypotheses. 
	
	We consider the class of differential systems
	\begin{align}
		\begin{split}	\label{pertorb-lineal-2}
			&\dot \theta = 1 + \varepsilon (a_0 + a_1 \theta + a_2 \varphi + a_3 r), \\
			& \dot \varphi = \varepsilon (b_0 + b_1 \theta + b_2 \varphi + b_3 r),  \\
			& \dot r = r- 1 + \varepsilon (c_0 + c_1 \theta + c_2 \varphi + c_3 r).
		\end{split}
	\end{align}
	where $a_i$, $b_i$ and $c_i$, for $i=0, \dots , 3$ are real numbers and with $\varepsilon >0$ being a small parameter. Note that this is the more general linear perturbation of system \eqref{sist-lineal-1}.
	For the class of systems \eqref{pertorb-lineal-2} we have the following result. 
	
	\begin{theorem} \label{th-pertorb-2} 
		For sufficiently small $\varepsilon >0$ the linear differential system \eqref{pertorb-lineal-2}  has a limit cycle bifurcating from a periodic orbit of system \eqref{sist-lineal-1} provided that $a_1b_2-a_2b_1 \neq 0$. Moreover this limit cycle bifurcates from the periodic orbit of system \eqref{sist-lineal-1} parameterized by $(\theta(t), \varphi(t), r(t)) = (\theta_0 + t, \varphi_0, 1)$, with 
		\begin{align*} 
			\theta_0 &= \frac{a_2(b_0+b_3+b_1 \pi)-b_2(a_0+a_3+a_1 \pi)}{a_1b_2-a_2b_1}, \\
			\varphi_0  &= \frac{b_1(a_0+a_3+a_1 \pi)-a_1(b_0+b_3+b_2 \pi)}{a_1b_2-a_2b_1}.
		\end{align*}
	\end{theorem}
	
	Theorem \ref{th-pertorb-2} is proved in section \ref{sec-dem-th1}. 
	
	We remark that the existence of the limit cycle for system \eqref{pertorb-lineal-2} does not depend on the perturbation of the $\dot r$ equation. 
	
	As an example of the previous result, consider  the system
	\begin{equation}  \label{pertorb-lineal-1}
		\dot \theta = 1 + \varepsilon a \varphi, \quad
		\dot \varphi = \varepsilon b \theta, \quad
		\dot r = r-1, 
	\end{equation}
	with $a, b \in \R$ and $\varepsilon >0$. In this case the sphere $\{ r=1 \}$ is still an invariant manifold. Appliying Theorem \ref{th-pertorb-2} with $a_2 = a, \, b_1 = b$ and the rest of the coefficients of the perturbation being zero,  we find that system \eqref{pertorb-lineal-1} has a limit cycle bifurcating form the periodic orbit of system \eqref{sist-lineal-1} parameterized by $(\theta (t), \varphi(t), r(t)) = (-\pi + t, 0, 1)$. That is, there is a limit cycle bifurcating from the periodic orbit corresponding to the equator of the sphere $\{r=1\}$ of system \eqref{sist-lineal-1}.  Moreover this limit cycle is still contained in the sphere $\{r=1\}$.
	
	Next we consider linear differential systems defined on higher dimensional manifolds.  We take $M = \mathbb{S}^2 \times \mathbb{S}^2 \times \R$ and 
	\begin{equation} 
		\label{sist-lineal-high}
		\dot \theta = 1, 
		\quad \dot \varphi = 0, 
		\quad \dot \nu = 1, 
		\quad \dot \phi = 0, 
		\quad \dot r = r-1,
	\end{equation}
	for  $(\theta, \varphi, \nu, \phi, r) \in M$, with $\theta, \nu \in [-\pi, \pi)$ and $\varphi, \phi \in (-\pi/2, \pi/2)$, 
	and with 
	$ \dot \theta = 0$ when $\varphi = \pm \pi / 2
	$
	and 
	$\dot \nu = 0$ when $\phi = \pm \pi / 2. $

	The general solution of system \eqref{sist-lineal-high} is
	$$
	\theta(t) = \theta_0 + t, \quad \varphi(t) = \varphi_0,  \quad \nu(t) = \nu_0 + t, \quad \phi(t) = \phi_0, \quad r(t) = (r_0 -1) e^t +1,
	$$
	and thus the product of spheres $\{r=1\} \cong (\mathbb{S}^2)^2 $ is an invariant manifold foliated by periodic orbits of period $2 \pi$, except for the four points $\{r=1, \varphi =\pm \pi/2 , \phi = \pm \pi/2 \}$,  which are equilibrium points. 
	
	We consider the most general perturbation of the differential system \eqref{sist-lineal-high} inside the class of all linear differential systems, namely
	\begin{align} 
		\begin{split} \label{pertorb-lineal-high}
			& \dot \theta = 1 + \varepsilon (a_0 +a_1 \theta + a_2 \varphi + a_3 \nu + a_4 \phi + a_5  r), \\
			& \dot \varphi = \varepsilon (b_0 +b_1 \theta + b_2 \varphi + b_3 \nu + b_4 \phi + b_5  r) \\
			& \dot \nu = 1 + \varepsilon(c_0 +c_1 \theta + c_2 \varphi + c_3 \nu + c_4 \phi + c_5  r), \\
			& \dot \phi = \varepsilon (d_0 +d_1 \theta + d_2 \varphi + d_3 \nu + d_4 \phi + d_5  r)  \\
			& \dot r = r-1 +  \varepsilon (e_0 +e_1 \theta + e_2 \varphi + e_3 \nu + e_4 \phi + e_5  r), 
		\end{split}
	\end{align}
	with $a_i, b_i, c_i, d_i, e_i \in \R$ for $i=0,\ldots,5$,
	and with $\varepsilon >0$ being a small parameter.  In the following result we give sufficient conditions on the coefficients of system \eqref{pertorb-lineal-high} in order that there is a limit cycle bifurcating from a periodic orbit of the corresponding unperturbed system.
	
	\begin{theorem} \label{th-pertorb-high} For sufficiently small $\varepsilon >0$ the differential system \eqref{pertorb-lineal-high}  has a limit cycle bifurcating from a periodic orbit of system \eqref{sist-lineal-high} provided that
		$$
		\det	\begin{pmatrix} 
			a_1 & a_2 & a_3 & a_4 \\
			b_1 & b_2 & b_3 & b_4 \\
			c_1 & c_2 & c_3 & c_4 \\
			d_1 & d_2 & d_3 & d_4 
		\end{pmatrix}  \neq 0.
		$$
		Moreover this limit cycle bifurcates from the periodic orbit of system \eqref{sist-lineal-high} parameterized by $(\theta(t), \varphi(t), \nu(t), \phi(t), r(t)) = (\theta_0 + t, \varphi_0,  \nu_0 + t, \phi_0, 1)$, where $(\theta_0, \,  \varphi_0, \,  \nu_0, \, \phi_0)$ is the unique solution of the linear system
		\begin{align*} 
			\begin{split} 
				& a_1 \theta_0 + a_2 \varphi_0 + a_3 \nu_0 + a_4 \phi_0 = -a_0 -a_1 \pi -a_3 \pi -a_5,  \\
				& b_1 \theta_0 + b_2 \varphi_0 + b_3 \nu_0 + b_4 \phi_0 = -b_0 -b_1 \pi -b_3 \pi -b_5,  \\
				& c_1 \theta_0 + c_2 \varphi_0 + c_3 \nu_0 + c_4 \phi_0 = -c_0 -c_1 \pi -c_3 \pi -c_5,  \\
				& d_1 \theta_0 + d_2 \varphi_0 + d_3 \nu_0 + d_4 \phi_0 = -d_0 -d_1 \pi -d_3 \pi -d_5. \\
			\end{split}
		\end{align*}
	\end{theorem}
	Theorem \ref{th-pertorb-high} is proved in section \ref{sec-dem-th2}. 
	
	Finally we consider the linear differential system defined in $M = \R^2 \times \mathbb{S}^2$, for $(x, y, \theta, \varphi) \in \R^2 \times \mathbb{S}^2$, with $\theta \in [-\pi, \pi)$ and $\varphi \in (-\pi/2, \pi/2)$, given by
	\begin{equation}  \label{sist-lineal-r2s2}
		\dot x = -y, 
		\quad \dot y = x, 
		\quad \dot \theta = 1, 
		\quad \dot \varphi = 0, 
	\end{equation}
	and with $\dot \theta = 0$ in the planes $P_1 =\{\varphi =-\pi/2\}$ and  $P_2 =\{\varphi =\pi/2\}$, which are invariant. 
	The general solution of  system \eqref{sist-lineal-r2s2} is 
	$$
	x(t) = x_0 \cos t - y_0 \sin t, \quad y(t) = x_0 \sin t + y_0 \cos t , \quad \theta(t) = \theta_0 + t, \quad \varphi(t) = \varphi_0,
	$$
	and therefore the whole phase space is filled by periodic orbits of period $2 \pi$, except for the two equilibrium points $(x, y, \theta, \varphi) = (0,0,\theta, -\pi/2)$ and $(x, y, \theta, \phi) = (0,0,\theta, \pi/2)$.
	
	We consider the most general linear perturbation of system \eqref{sist-lineal-r2s2} and we study the existence of limit cycles bifurcating  from the periodic orbits of system \eqref{sist-lineal-r2s2}.
	
	Let
	\begin{align} 
		\begin{split} \label{pertorb-lineal-r2s2}
			& \dot x = -y + \varepsilon (a_0 + a_1 x+ a_2 y+ a_3 \theta + a_4 \varphi), \\
			& \dot y = x + \varepsilon (b_0 + b_1 x+ b_2 y+ b_3 \theta + b_4 \varphi), \\
			& \dot \theta = 1 + \varepsilon (c_0 + c_1 x+ c_2 y+ c_3 \theta + c_4 \varphi), \\
			& \dot \varphi =  \varepsilon (d_0 + d_1 x+ d_2 y+ d_3 \theta + d_4 \varphi), \\
		\end{split}
	\end{align}
	be the perturbed system,
	with $a_i, b_i, c_i, d_i\in \R$ for $i=0,\ldots,4$, and where $\varepsilon >0$ is a small parameter. For this linear differential system  we  have  the following result.

	\begin{theorem} \label{th-pertorb-r2s2}For sufficiently small $\varepsilon >0$ the linear differential system \eqref{pertorb-lineal-r2s2}  has a limit cycle bifurcating from a periodic orbit of system \eqref{sist-lineal-r2s2} provided that 
		$$
		\det \begin{pmatrix}
			b_2+a_1 & a_2-b_1 \\
			b_1-a_2 & b_2+a_1
		\end{pmatrix} \neq 0 \quad \text{and} \quad
		\det \begin{pmatrix}
			c_3 & c_4 \\
			d_3 & d_4
		\end{pmatrix} \neq 0.
		$$	
		Moreover this limit cycle bifurcates form the periodic orbit of system \eqref{sist-lineal-r2s2} passing through the point $(x_0, y_0, \theta_0, \varphi_0)$ where
		\begin{align*}
			x_0 & = \frac{(2 b_2+2 a_1) b_3-2a_3 b_1+2 a_2a_3}{b_2^2 +b_1^2 + a_2^2+a_1^2 +2 a_1 b_2-2a_2b_1}, \\
			y_0 &= - \frac{(2b_1-2a_2)b_3+2a_3 b_2+2a_1a_3}{b_2^2 +b_1^2 + a_2^2+a_1^2 +2 a_1 b_2-2a_2b_1}, \\
			\theta_0 &= - \frac{(\pi c_3+c_0)d_4-\pi c_4d_3-c_4 d_0}{c_3 d_4 - c_4 d_3}, \\
			\varphi_0 &=  \frac{c_0 d_3 - c_3 d_0}{c_3 d_4 - c_4 d_3}.
		\end{align*}
	\end{theorem}
	Theorem \ref{th-pertorb-r2s2} is proved in section \ref{sec-dem-th3}. 
	
	As an example consider the system
	\begin{equation} 
		\label{pertorb-lineal-3}
		\dot x = -y + \varepsilon a y, \quad
		\dot y = x + \varepsilon b x, \quad
		\dot \theta = 1 + \varepsilon c \varphi, \quad
		\dot \varphi =  \varepsilon d \theta, 
	\end{equation}
	with  $a, b, c , d \in \R$, and $\varepsilon >0$. Appliying Theorem \ref{th-pertorb-r2s2} with $a_2 = a, \, b_1 = b, \, c_4 =c, \, d_3 = d$ and the rest of the coefficients of the perturbation being zero,  we obtain that system \eqref{pertorb-lineal-3} has a limit cycle bifurcating form the periodic orbit of system \eqref{sist-lineal-r2s2}  passing through the point $(x_0, y_0, \theta_0, \varphi_0) = (0,0, -\pi, 0)$,  provided that $(a-b)cd \neq 0$. That is, here the limit cycle bifurcates from the periodic orbit corresponding to the equator of the invariant sphere $\{x=y=0\}$ of system \eqref{sist-lineal-r2s2}. 
	
	The key tool that we use for proving Theorems \ref{th-pertorb-2} - \ref{th-pertorb-r2s2} is the averaging theory. For a general introduction to this theory, see the books \cite{san-ver-mur, verhulst}. As one can see in the proofs of Theorems \ref{th-pertorb-2} - \ref{th-pertorb-r2s2}, our method based on the averaging theory can produce at most one limit cycle for the studied systems. Therefore the following open question is natural.
	
	\textbf{Open question.} Let $m$ and $n$ be two non-negative integers.  Is it true that a linear  vector field on the manifold  $(\mathbb{S}^m)^m \times \R^n$ can have at most one limit cycle?
	
	A similar open question was stated in \cite{lli-zhang} concerning linear vector fields on the manifold $(\mathbb{S}^1)^m \times \R^n$.

	\section{Basic results on the averaging theory} \label{sec-teoria-limitcycles}
	
	In this section we state some basic results from the averaging theory that will be used to prove the main results of the paper.
	
	Let $M$ be a smooth connected manifold of dimension $n$, and let $F_0, \,  F_1: \R \times M \to \R^n$ and $F_2  : \R \times M \times [0,  \varepsilon_0) \to \R^n$ be $C^2$, $T$-periodic functions.
	Given the differential system 
	\begin{equation} \label{sist-fonamental-gral}
		\dot x (t) = F_0(t,x),
	\end{equation}
	we consider a perturbation of this system of the form
	\begin{equation} \label{sist-pertorbat-gral}
		\dot x (t) = F_0(t,x) + \varepsilon F_1(t,x) + \varepsilon^2 F_2 (t, x, \varepsilon).
	\end{equation} 
	
	The objective is to study the bifurcation of $T$-periodic solutions of system  \eqref{sist-pertorbat-gral} for $\varepsilon>0$ small enough. A solution to this problem is given by the averaging theory.

	We assume that there exists $k \leq n$ such that $M = M_k \times M_{n-k}$, where $M_k$ is a manifold of dimension $k$ and $M_{n-k}$ is a manifold of dimension $n-k$, and that
	the unperturbed system, namely system \eqref{sist-fonamental-gral}, contains an  open set, $V \subseteq M_k$, such that $\overline V$ is filled with periodic solutions all of them with the same period. 
	Such a set is called \emph{isochronous}.

	Let $x(t, z, \varepsilon)$ be the solution of system \eqref{sist-pertorbat-gral} such that $x(0, z, \varepsilon) =z$. We write the linearization of the unperturbed system \eqref{sist-fonamental-gral} along the solution $x(t,z,0)$ as
	\begin{equation} \label{variacional-gral}
		\dot y =  D_x F_0 (t,  x(t, z, 0)) y,
	\end{equation}
	and we denote by $\mathcal{M}_z(t)$ the fundamental matrix of the linear differential system \eqref{variacional-gral} such that $\mathcal{M}_z(0)$ is the $n \times n$ identity matrix, and by $\xi : M = M_k \times M_{n-k} \to M_k$ the projection of $M$ onto its first $k$ coordinates, that is, $\xi(x_1, \dots, x_n) = (x_1, \dots, x_k)$.
	
	The following results give sufficient conditions for the existence of limit cycles for a system of the form \eqref{sist-pertorbat-gral} bifurcating from the periodic orbits of system \eqref{sist-fonamental-gral}.

	\begin{theorem} \label{th-averaging-gral}
		Let $V \subseteq M_k$ be an open and bounded set, and let $\beta_0: \overline{V} \to M_{n-k}$ be a $C^ 2$ function. Assume
		
		(i) $\mathcal{Z}=  \{ z_\alpha =(\alpha, \beta_0(\alpha)) \,  : \,  \alpha \in  \overline{V}  \} \subset M$ and  for each  $z_\alpha \in \mathcal{Z}$ the solution $x(t, z_\alpha,  0)$ of system \eqref{sist-fonamental-gral} is T-periodic.
		
		(ii) For each $z_\alpha \in \mathcal{Z}$, there is a  fundamental matrix $\mathcal{M}_{z_\alpha}(t)$ of system \eqref{variacional-gral}  such that the matrix $\mathcal{M}_{z_\alpha}^{-1}(0) - \mathcal{M}_{z_\alpha}^{-1}(T)$ has the $k \times (n-k)$ zero matrix in the upper right  corner, and  a  $(n-k) \times (n-k)$ matrix $\Delta_\alpha$ in the lower right corner with $\det (\Delta_\alpha) \neq 0$.
		
		Consider the function $\mathcal{F}: \overline V \to \R^k$ defined by 
		\begin{equation*} \label{integral-averaging}
			\mathcal{F}(\alpha) = \xi \bigg{(}  \int_0^T \mathcal{M}_{z_\alpha}^{-1}(t) F_1(t, x(t, z_\alpha, 0)) \,  dt \bigg{)}.
		\end{equation*}
		
		If there exists $a \in V$ with $\mathcal{F}(a)=0$ and with $\det  (\mathcal{DF}(a)) \neq 0$, then there is a limit  cycle $x(t, \varepsilon)$ of period $T$ of system \eqref{sist-pertorbat-gral}  such that $x (0, \varepsilon) \to z_a$ as $\varepsilon \to 0$.
	\end{theorem}
	
	The result given by Theorem \ref{th-averaging-gral} can be found in the books of Malkin \cite{malkin} and Rosseau \cite{rosseau}. For a shorter proof, see \cite{buica-francoise-llibre}. There the result is proved in $\R^n$, but it can be easily extended to a manifold $M$.

	The next result allows to determine the existence of limit cycles in a system of the form \eqref{sist-pertorbat-gral} in the case when there exists an open set, $V \subset M$, such that for all $z \in \overline V$, the solution $x(t, z,0)$  is $T$-periodic. 
	
	\begin{theorem} \label{th-averaging-2}
		Let $V \subseteq M$ be an open and bounded set with $\overline V \subseteq M$, and assume that for all $z \in \overline V$ the solution $x(t, z, 0)$ of system \eqref{sist-pertorbat-gral} is $T$-periodic. Consider the function $\mathcal{F}: \overline V \to \R^n$ defined by 
		\begin{equation*} \label{integral-averaging-particular}
			\mathcal{F}(z) =   \int_0^T M_{z}^{-1}(t) F_1(t, x(t, z, 0)) \,  dt.
		\end{equation*}	
		If there exists $a \in V$ with $\mathcal{F}(a)=0$ and with $\det  (\mathcal{DF}(a)) \neq 0$, then there is a limit  cycle $x(t, \varepsilon)$ of period $T$ of system \eqref{sist-pertorbat-gral}  such that $x (0, \varepsilon) \to a$ as $\varepsilon \to 0$.
	\end{theorem}
	
	For the proof of Theorem \ref{th-averaging-2} see Corollary $1$ of \cite{buica-francoise-llibre}.
	
	\section{Proof of Theorem \ref{th-pertorb-2}} \label{sec-dem-th1}
	
	We use the result from averaging theory given in Theorem \ref{th-averaging-gral} to deduce the existence of a limit cycle of system \eqref{pertorb-lineal-2}, for some $\varepsilon>0$ small enough, bifurcating from a periodic orbit of the same system with $\varepsilon=0$.
	
	Since the general solution of the differential system \eqref{sist-lineal-1}, corresponding to system \eqref{pertorb-lineal-2} with $\varepsilon = 0$, is given by 
	$$
	\theta(t) = \theta_0 + t, \quad \varphi(t) = \varphi_0, \quad r(t) = (r_0 -1) e^t +1, 
	$$
	it is clear that all the periodic solutions of that system are parameterized by 
	$$
	\theta(t) = \theta_0 + t, \quad \varphi(t) = \varphi_0, \quad r(t) = 1,
	$$
	with $(\theta_0, \varphi_0) \in \mathbb{S}^2 \backslash \{\varphi_0 = \pm \pi/2\}$. Then, all the periodic solutions have period $2 \pi$ and they fill the invariant sphere $\{r=1\}$ except for the poles, which are equilibrium points.
	
	Therefore, for applying Theorem \ref{th-averaging-gral} we take $M = \mathbb{S}^2 \times \R$ and
	\begin{align}
		\begin{split} \label{cond-hipotesis-1}
			k & =2, \, n  = 3, \\
			M_k &= M_2 = \{(\theta, \varphi, r) \in  M \, : \, r=1 \} \cong  \mathbb{S}^2, \\
			x & = (\theta, \varphi, r), \\
			\alpha & = (\theta_0, \varphi_0), \\
			\beta_0(\alpha) & = \beta_0(\theta_0, \varphi_0) = 1, \\
			z_\alpha & = (\alpha, \beta_0(\alpha)) = (\theta_0, \varphi_0, 1), \\
			V  &= \{  (\theta, \varphi, r) \in M \, : \, r=1, \,  \varphi \in (-\tfrac{\pi}{2} + \delta_0, \tfrac{\pi}{2} - \delta_0 ) \} \\
			& \quad  \text{ with $\delta_0 >0$ small enough such that } \\
			&  \quad \,  \varphi^* := \frac{b_1(a_0+a_3+a_1 \pi)-a_1(b_0+b_3+b_2 \pi)}{a_1b_2-a_2b_1} \in (-\tfrac{\pi}{2} + \delta_0, \tfrac{\pi}{2} - \delta_0 ), \\
			\mathcal{Z} & = \overline V \times \{r=1\}, \\
			x(t, z_\alpha, 0) &= (\theta_0 + t, \varphi_0, 1), \\
			F_0(t, x) & = (1, 0, r-1), \\
			F_1(t, x) & =  (a_0 + a_1 \theta + a_2 \varphi + a_3 r, \, b_0 + b_1 \theta + b_2 \varphi + b_3 r, \, c_0 + c_1 \theta + c_2 \varphi + c_3 r ), \\
			F_2(t, x, \varepsilon) & = 0, \\
			T & = 2 \pi,
		\end{split}
	\end{align}
	where we took $V \subset  M_2$ as an open subset that contains the periodic orbit for which it bifurcates a limit cycle, as we shall see next. 
	
	The fundamental matrix $\mathcal{M}_{z_\alpha}(t)$ with $\mathcal{M}_{z_\alpha}(0) = Id$ of system \eqref{variacional-gral} with $F_0$ and $x(t, z_\alpha, 0)$ described above is the matrix $\mathcal{M}_{z_\alpha}(t) = \exp(D_xF_0 \, t)$, \emph{i.e.}
	\begin{equation*} \label{fundamental-pertorb-1}
		\mathcal{M}_{z_\alpha}(t) = \begin{pmatrix}
			1  & 0 & 0 \\
			0 & 1 & 0 \\
			0 & 0 & e^t
		\end{pmatrix}.
	\end{equation*}
	Note that since $F_0$ defines a linear differential system, the fundamental matrix $\mathcal{M}_{z_\alpha}(t)$ is independent of the initial conditions $z_\alpha$. 
	
	We also have
	$$
	\mathcal{M}_{z_\alpha}^{-1} (0) - \mathcal{M}_{z_\alpha}^{-1} (2 \pi) = \begin{pmatrix}
		0 &0 &0 \\
		0 & 0 & 0 \\
		0 & 0 & 1-e^{-2\pi}
	\end{pmatrix},
	$$
	and therefore, all the assumptions in the in the statement of Theorem \ref{th-averaging-gral} are satisfied.

	With the described setting, the function $\mathcal{F}(\alpha) = \mathcal{F}(\theta_0, \varphi_0)$ from the statement of Theorem \ref{th-averaging-gral} associated with system \eqref{pertorb-lineal-2} is 
	\begin{align*}
		\mathcal{F}(\theta_0, \varphi_0) &= \xi \bigg{(}  \int_0^{2 \pi} \mathcal{M}_{z_\alpha}^{-1}(t) F_1(\theta_0 + t, \varphi_0, 1) \,  dt \bigg{)}  \\
		& = 2 \pi (a_0 + a_1(\theta_0+\pi) + a_2 \varphi_0 + b_3, \, b_0 + b_1(\theta_0+\pi) + b_2 \varphi_0 + b_3). 
	\end{align*}
	We have $\det (D\mathcal{F}) = 4 \pi^ 2 (a_1b_2-a_2b_1)$, and therefore $\det (D\mathcal{F}) \neq 0$ for all $(\theta_0, \varphi_0) \in V$. Thus, the only solution of $\mathcal{F} = 0$ is given by 
	\begin{align} 
		\begin{split} \label{theta-phi-2}
			\theta_0 &= \frac{a_2(b_0+b_3+b_1 \pi)-b_2(a_0+a_3+a_1 \pi)}{a_1b_2-a_2b_1}, \\
			\varphi_0  &= \frac{b_1(a_0+a_3+a_1 \pi)-a_1(b_0+b_3+b_2 \pi)}{a_1b_2-a_2b_1}.
		\end{split}
	\end{align}
	Note that such solution $(\theta_0, \varphi_0)$, where $\varphi_0 = \varphi^*$, is contained in the set $V$ described in \eqref{cond-hipotesis-1}.
	
	Hence, by Theorem \ref{th-averaging-gral}, if $\varepsilon >0$ is small enough, there is a periodic solution, $(\theta(t, \varepsilon), \varphi(t, \varepsilon), r(t, \varepsilon))$, of system \eqref{pertorb-lineal-1}, which is a limit cycle, and such that 
	$$
	(\theta(0, \varepsilon), \varphi(0, \varepsilon), r(0, \varepsilon)) \to (\theta_0, \varphi_0, 1),
	$$ 
	
	when $\varepsilon \to 0$, and where $\theta_0$ and $\varphi_0$ are given in \eqref{theta-phi-2}.

	\section{Proof of Theorem \ref{th-pertorb-high}} \label{sec-dem-th2}
	
	We use the result from averaging theory given in Theorem \ref{th-averaging-gral} to prove that, for some $\varepsilon>0$ small enough, there exist a limit cycle of system \eqref{pertorb-lineal-high} bifurcating from a periodic orbit of the same system with $\varepsilon=0$.
	
	Since the general solution of system \eqref{pertorb-lineal-high} with $\varepsilon = 0$ (that is, the one of system \eqref{sist-lineal-high}), is 
	$$
	\theta(t) = \theta_0 + t, \quad \varphi(t) = \varphi_0, \quad \nu(t) = \nu_0 + t, \quad \phi(t) = \phi_0, \quad r(t) = (r_0 -1) e^t +1, 
	$$
	then all the periodic solutions of that system are  
	$$
	\theta(t) = \theta_0 + t, \quad \varphi(t) = \varphi_0,  \quad \nu(t) = \nu_0 + t, \quad \phi(t) = \phi_0,  \quad r(t) = 1,
	$$
	with $(\theta_0, \varphi_0, \nu_0, \phi_0) \in \mathbb{S}^2 \backslash \{\varphi_0 = \pm \pi/2\} \times \mathbb{S}^2 \backslash \{\varphi_0 = \pm \pi/2\}$. That is, the periodic solutions fill the invariant manifold $\{r=1\}$ except for the four equilibrium points $\{\varphi= \pm \pi/2, \, \phi = \pm \pi/2 \}$, and they have all period $2 \pi$. 
	
	For applying Theorem \ref{th-averaging-gral} we take $M = (\mathbb{S}^2)^2 \times \R$ and 
	\begin{align}
		\begin{split} \label{cond-hipotesis-high}
			k & =4, \,	n  = 5, \\
			M_k & = M_4 = \{ \theta, \varphi,  \nu, \phi, r \in M \, : \, r=1 \} \cong (\mathbb{S}^2)^2, \\
			x & = (\theta, \varphi, \nu, \phi, r), \\
			\alpha & = (\theta_0, \varphi_0, \nu_0, \phi_0), \\
			\beta_0(\alpha) & = \beta_0(\theta_0, \varphi_0, \nu_0, \phi_0) = 1, \\
			z_\alpha & = (\alpha, \beta_0(\alpha)) = (\theta_0, \varphi_0, \nu_0, \phi_0, 1), \\
			V &= \{  (\theta, \varphi, \nu, \phi, r) \in M  \,  : \, r=1, \, \varphi \in (-\tfrac{\pi}{2}+ \delta_0, \tfrac{\pi}{2}-\delta_0)  \}  \\
			& \quad  \text{ with $\delta_0 >0$ small enough such that $\varphi_0, \, \phi_0$ satisfying \eqref{sist-lin-cycle} satisfy} \\
			&  \quad \, \varphi_0, \, \phi_0 \in (-\tfrac{\pi}{2} + \delta_0, \tfrac{\pi}{2} - \delta_0 ), \\
			\mathcal{Z} & = \overline V \times \{ r=1\}, \\
			x(t, z_\alpha, 0) &= (\theta_0 + t, \varphi_0, \nu_0 +t, \phi_0, 1), \\
			F_0(t, x) & = (1, 0, 1, 0, r-1), \\
			F_1(t, x) & = \begin{pmatrix}
				a_0 +a_1 \theta + a_2 \varphi + a_3 \nu + a_4 \phi + a_5  r \\
				b_0 +b_1 \theta + b_2 \varphi + b_3 \nu + b_4 \phi + b_5  r \\
				c_0 +c_1 \theta + c_2 \varphi + c_3 \nu + c_4 \phi + c_5  r\\
				d_0 +d_1 \theta + d_2 \varphi + d_3 \nu + d_4 \phi + d_5  r \\
				e_0 +e_1 \theta + e_2 \varphi + e_3 \nu + e_4 \phi + e_5  r
			\end{pmatrix}, \\
			F_2(t, x, \varepsilon) & = 0, \\
			T & = 2 \pi,
		\end{split}
	\end{align}
	where we chose $V \subset  M_4$ as an open subset that contains the periodic orbit for which it bifurcates a limit cycle, as we shall see next. 
	
	The fundamental matrix $\mathcal{M}_{z_\alpha}(t)$ with $\mathcal{M}_{z_\alpha}(0)=Id$, of system \eqref{variacional-gral} with $F_0$ and $x(t, z_\alpha, 0)$ described above is the matrix $\mathcal{M}_{z_\alpha}(t) = \exp(D_xF_0 \, t)$, \emph{i.e.},
	\begin{equation*} \label{fundamental-pertorb-2}
		\mathcal{M}_{z_\alpha}(t) = \begin{pmatrix}
			1&0&0&0&0 \\
			0&1&0&0&0 \\
			0&0&1&0&0 \\
			0&0&0&1&0 \\
			0 & 0 & 0&0 & e^t
		\end{pmatrix}.
	\end{equation*}
	
	We also have
	$$
	\mathcal{M}_{z_\alpha}^{-1} (0) - \mathcal{M}_{z_\alpha}^{-1} (2 \pi) = \begin{pmatrix}
		0 &0 &0& 0 &0 \\
		0 &0 &0& 0 &0 \\
		0 &0 &0& 0 &0 \\
		0 &0 &0& 0 &0 \\
		0 & 0 & 0 & 0 &1-e^{-2\pi}
	\end{pmatrix},
	$$
	and therefore, all the assumptions in the statement of Theorem \ref{th-averaging-gral} are satisfied. 
	
	With the described setting, the function $\mathcal{F}(\alpha) = \mathcal{F}(\theta_0, \varphi_0, \nu_0, \phi_0)$ in the statement of Theorem \ref{th-averaging-gral} associated with system \eqref{pertorb-lineal-high} is 
	\begin{align*}
		\mathcal{F}(\theta_0, \varphi_0, \nu_0, \phi_0) & = \xi \bigg{(}  \int_0^{2 \pi} \mathcal{M}_{z_\alpha}^{-1}(t) F_1(\theta_0 + t, \varphi_0, \nu_0 + t, \phi_0, 1) \,  dt \bigg{)}  = (\mathcal{F}_1, \mathcal{F}_2, \mathcal{F}_3, \mathcal{F}_4 ),
	\end{align*}
	with 
	\begin{align*}
		\mathcal{F}_1&= 2 \pi (a_0 + a_1 \theta_0 +  a_1  \pi + a_2 \varphi_0 + a_3 \nu_0 + a_3 \pi + a_4 \phi_0 + a_5), \\
		\mathcal{F}_2&= 2 \pi (b_0 + b_1 \theta_0 +  b_1  \pi + b_2 \varphi_0 + b_3 \nu_0 + b_3 \pi + b_4 \phi_0 + b_5), \\
		\mathcal{F}_3&= 2 \pi (c_0 + c_1 \theta_0 +  c_1  \pi + c_2 \varphi_0 + c_3 \nu_0 + c_3 \pi + c_4 \phi_0 + c_5) ,\\
		\mathcal{F}_4&= 2 \pi (d_0 + d_1 \theta_0 +  d_1  \pi + d_2 \varphi_0 + d_3 \nu_0 + d_3 \pi + d_4 \phi_0 + d_5) .
	\end{align*}
	Also, we have 
	$$
	\det (D\mathcal{F}) = 16 \pi^4 \,  \det \begin{pmatrix} 
		a_1 & a_2 & a_3 & a_4 \\
		b_1 & b_2 & b_3 & b_4 \\
		c_1 & c_2 & c_3 & c_4 \\
		d_1 & d_2 & d_3 & d_4 
	\end{pmatrix}  \neq 0,
	$$
	by assumption. The initial conditions $(\theta_0, \varphi_0, \nu_0, \phi_0)$ such that $\mathcal{F} (\theta_0, \varphi_0, \nu_0, \phi_0) =0$ are the solutions of the linear system 
	\begin{align} 
		\begin{split} \label{sist-lin-cycle}
			& a_1 \theta_0 + a_2 \varphi_0 + a_3 \nu_0 + a_4 \phi_0 = -a_0 -a_1 \pi -a_3 \pi -a_5,  \\
			& b_1 \theta_0 + b_2 \varphi_0 + b_3 \nu_0 + b_4 \phi_0 = -b_0 -b_1 \pi -b_3 \pi -b_5,  \\
			& c_1 \theta_0 + c_2 \varphi_0 + c_3 \nu_0 + c_4 \phi_0 = -c_0 -c_1 \pi -c_3 \pi -c_5,  \\
			& d_1 \theta_0 + d_2 \varphi_0 + d_3 \nu_0 + d_4 \phi_0 = -d_0 -d_1 \pi -d_3 \pi -d_5. \\
		\end{split}
	\end{align}
	
	Since $	\det (D\mathcal{F})  \neq 0$,  system \eqref{sist-lin-cycle} has a unique solution. Note that such solution $(\theta_0, \varphi_0, \nu_0, \phi_0)$ is contained in the set $V$ described in \eqref{cond-hipotesis-high}.
	
	Hence, by Theorem \ref{th-averaging-gral}, if $\varepsilon >0$ is small enough, there is a periodic solution, 
	$$
	(\theta(t, \varepsilon), \varphi(t, \varepsilon), \nu(t, \varepsilon), \phi(t, \varepsilon), r(t, \varepsilon)),
	$$ of system \eqref{pertorb-lineal-high}, which is a limit cycle, and such that 
	$$
	(\theta(0, \varepsilon), \varphi(0, \varepsilon),  \nu(0, \varepsilon), \phi(0, \varepsilon),  r(0, \varepsilon)) \to (\theta_0, \varphi_0, \nu_0, \phi_0, 1),
	$$ 
	when $\varepsilon \to 0$, and where $\theta_0, \, \varphi_0, \, \nu_0, $ and $ \phi_0$ are given by the unique solution of system \eqref{sist-lin-cycle}.

	\section{Proof of Theorem \ref{th-pertorb-r2s2}} \label{sec-dem-th3}
	
	Since the general solution of system \eqref{pertorb-lineal-r2s2} with $\varepsilon = 0$ is given by 
	\begin{equation*} \label{sol-th-r2s2}
		x(t) = x_0 \cos(t) - y_0 \sin(t), \quad y(t) = x_0 \sin(t) + y_0 \cos(t) , \quad \theta(t) = \theta_0 + t, \quad \varphi(t) = \varphi_0,
	\end{equation*}
	the whole phase space is filled by periodic solutions, except form the equilibrium points $(x, y, \theta, \varphi) = (0,0,\theta, -\pi/2)$ and $(x, y, \theta, \varphi) = (0,0,\theta, \pi/2)$. 
	Hence, the periodic solutions of the differential system \eqref{sist-lineal-r2s2} fill an open set of the phase space $M = \R^2 \times \mathbb{S}^2$.
	
	To prove Theorem \ref{th-pertorb-r2s2} we use the result given in Theorem \ref{th-averaging-2} to deduce that there exist a limit cycle of system \eqref{pertorb-lineal-r2s2}, for some $\varepsilon>0$ small enough, bifurcating from the periodic orbits of the same system with $\varepsilon=0$.
	
	To clarify the notation, here the solution $x(t, z, 0)$ from the statement of Theorem \ref{th-averaging-2} will be denoted by $\mathbf{x}(t, z, 0)$, and $x$ will denote the first variable in the phase space.

	For applying Theorem \ref{th-averaging-2} we take $M= \R^2 \times \mathbb{S}^2$ and 
	\begin{align}
		\begin{split} \label{cond-hipotesis-2}
			\mathbf{x} & = (x, y, \theta, \varphi), \\
			z & = (x_0, y_0, \theta_0, \varphi_0), \\
			\mathbf{x}(t, z, 0) &= (x(t), y(t), \theta(t), \varphi(t)) \text{ given by \eqref{sol-th-r2s2}} \\
			F_0(t, x) & = (-y,x,1,0), \\
			F_1(t, x) & = \begin{pmatrix}
				a_0 +a_1 x + a_2 y + a_3 \theta + a_4 \varphi  \\
				b_0 +b_1 x + b_2 y + b_3 \theta + b_4 \varphi  \\
				c_0 +c_1 x + c_2 y + c_3 \theta + c_4 \varphi  \\
				d_0 +d_1 x + d_2 y + d_3 \theta + d_4 \varphi  \\
			\end{pmatrix}, \\
			F_2(t, x, \varepsilon) & = 0, \\
			T & = 2 \pi, \\
			V  &= \{(x, y, \theta, \varphi) \in M \, :  \, \|(x, y)\|<1 + \kappa, \, \varphi \in (-\tfrac{\pi}{2}+ \delta_0, \tfrac{\pi}{2} - \delta_0) \}, \\
			& \quad \ \text{with } \ \kappa = \frac{2 \,  \sqrt{a_3^2+b_3^2}}{\sqrt{a_1^2+a_2^2+b_1^2+b_2^2 +  2 a_1b_2 - 2a_2b_1}}, \\
			& \quad  \text{ and with $\delta_0 >0$ small enough such that } \\
			&  \quad \,  \varphi^* := \frac{c_0 d_3 - c_3 d_0}{c_3 d_4 - c_4 d_3}\in (-\tfrac{\pi}{2} + \delta_0, \tfrac{\pi}{2} - \delta_0 ).
		\end{split}
	\end{align}
	where we chose $V \subset  M$ as an open subset that contains the periodic orbit for which it bifurcates a limit cycle, as we shall see next. 
	
	The fundamental matrix $\mathcal{M}_{z}(t)$ of system \eqref{variacional-gral} with $\mathcal{M}_z(0)=Id$ and with $F_0$ and $\mathbf{x}(t, z, 0)$ described in \eqref{cond-hipotesis-2} is given by 
	\begin{equation*} \label{fundamental-pertorb-r2s2}
		M_{z}(t) = \begin{pmatrix}
			\cos (t) & -\sin(t) &0 &0 \\
			\sin(t) & \cos (t) &0 &0 \\
			0&0&1&0 \\
			0&0&0&1
		\end{pmatrix}.
	\end{equation*}
	Therefore all the assumptions in the statement of Theorem \ref{th-averaging-2} are satisfied. 
	
	With the described setting the function $\mathcal{F}(z) = \mathcal{F}(x_0, y_0, \theta_0, \varphi_0)$ in the statement of Theorem \ref{th-averaging-2} associated with system \eqref{pertorb-lineal-r2s2}, namely,
	\begin{align*}
		\mathcal{F}(x_0, y_0, \theta_0, \varphi_0)& =   \int_0^{2 \pi} M_{z}^{-1}(t) F_1(t, \mathbf{x}(t,x,0)) \,  dt,
	\end{align*}
	is given by $\mathcal{F}=(\mathcal{F}_1, \mathcal{F}_2, \mathcal{F}_3, \mathcal{F}_4)$, which  after some straightforward computations can be written as 
	\begin{align*}
		\mathcal{F}_1 & =(\pi a_2-\pi b_1) y_0+(\pi b_2+ \pi a_1)x_0-2 \pi b_3, \\
		\mathcal{F}_2 & = (\pi b_2+\pi a_1)y_0+(\pi b_1- \pi a_2)x_0+2 \pi a_3, \\
		\mathcal{F}_3 &= 2 \pi c_3 \theta_0+2 \pi c_4 \varphi_0+2 \pi^2 c_3+2 \pi c_0, \\
		\mathcal{F}_4 &= 2 \pi d_3 \theta_0+2 \pi d_4 \varphi_0 +2 \pi^2 d_3+2 \pi d_0.
	\end{align*}
	
	Assuming that
	\begin{equation} \label{det-3-r2s2}
		\det (\mathcal{DF}) =
		\det \begin{pmatrix}
			\pi (b_2 + a_1) & \pi (a_2-b_1) & 0 &0 \\
			\pi (b_1 - a_2) & \pi (b_2+a_1) & 0 &0 \\
			0 & 0 & 2 \pi c_3 & 2 \pi c_4 \\
			0 & 0 & 2 \pi d_3 & 2 \pi d_4 \\
		\end{pmatrix} \neq 0,
	\end{equation}
	the linear system $(\mathcal{F}_1, \mathcal{F}_2, \mathcal{F}_3, \mathcal{F}_4) = (0,0,0,0)$ has a unique solution, given by
	\begin{align}
		\begin{split} \label{sol-punts-r2s2}
			x_0 & = \frac{(2 b_2+2 a_1) b_3-2a_3 b_1+2 a_2a_3}{b_2^2 +b_1^2 + a_2^2+a_1^2 +2 a_1 b_2-2a_2b_1}, \\
			y_0 &= - \frac{(2b_1-2a_2)b_3+2a_3 b_2+2a_1a_3}{b_2^2 +b_1^2 + a_2^2+a_1^2 +2 a_1 b_2-2a_2b_1}, \\
			\theta_0 &= - \frac{(\pi c_3+c_0)d_4-\pi c_4d_3-c_4 d_0}{c_3 d_4 - c_4 d_3}, \\
			\varphi_0 &=  \frac{c_0 d_3 - c_3 d_0}{c_3 d_4 - c_4 d_3}.
		\end{split}
	\end{align}
	Note that such solution $(x_0, y_0, \theta_0, \varphi_0)$, where $\varphi_0 = \varphi^*$, is contained in the set $V$ described in \eqref{cond-hipotesis-2}.
	
	The condition  \eqref{det-3-r2s2} 
	is clearly satisfied for all $(x_0, y_0, \theta_0, \varphi_0) \in V$  taking into account the assumptions in the statement of Theorem \ref{th-pertorb-r2s2}.
	
	Hence, by Theorem \ref{th-averaging-2}, there is a periodic solution $(x(t, \varepsilon), y(t, \varepsilon), \theta(t, \varepsilon), \varphi(t, \varepsilon))$ of system \eqref{pertorb-lineal-r2s2}, which is a limit cycle, and such that 
	$$
	(x(0, \varepsilon), y(0, \varepsilon), \theta(0, \varepsilon), \varphi(0, \varepsilon), r(0, \varepsilon)) \to (x_0, y_0,\theta_0, \varphi_0) 
	$$ 
	when $\varepsilon \to 0$, and where $x_0, \, y_0,\, \theta_0$ and $ \varphi_0$ are given in \eqref{sol-punts-r2s2}.
	
	\vspace{0.4cm}
	
	\textbf{\large{Acknowledgements.}} This work is supported by the Ministerio de Ciencia, Innovaci\'on y Universidades, Agencia Estatal de Investigaci\'on grants PID2019-104658GB-I00 and BES-2017-081570,  and the H2020 European Research Council grant MSCA-RISE-2017-777911. 
	
	\pagebreak


\begin{thebibliography}{99}
		
		\bibitem{arnold} V. I. Arnold, Ordinary Differential Equations, \emph{2nd printing of the 1992 edn} (2006). Berlin, Springer, Universitext. 
		
		\bibitem{buica-francoise-llibre} A. J. Buica, J. P. Françoise and J. Llibre, Periodic solutions of nonlinear periodic differential systems with a small parameter, \emph{Commun. Pure Appl. Anal} \textbf{6} (2007), 103-111.
		
		\bibitem{ChLlo} C. Christopher and N. G. Lloyd, Small-amplitude limit cycles in polynomial Liénard systems, \emph{Nonlinear Differ. Equ. Appl} \textbf{3} (1996), 183-190.
		
		\bibitem{fer-lli-teix} A. Ferragut, J. Llibre and M. A. Teixeira, Hyperbolic periodic orbits from the bifurcation of a four-dimensional nonlinear center, \emph{Int. J. Bifurcation Chaos Appl. Sci. Eng.} \textbf{17} (2007), 2159-2167.
		
		\bibitem{GiLliVia} H. Giacomini, J. Llibre and M. Viano, On the nonexistence, existence and uniqueness of limit cycles, \emph{Nonlinearity} \textbf{9} (1996), 501-516.
		
		\bibitem{HanLi} M. Han and L. Li, Lower bounds for the Hilbert number of polynomial systems, \emph{J. Differential Equations} \textbf{252} (2012), 3278-3304.
		
		\bibitem{hartman} P. Hartman, A lemma in the theory of structural stability of differential equations, \emph{Proc. Am. Math. Soc} \textbf{11} (1960), 610-620.
		
		\bibitem{Ilya} Yu. Ilyashenko, Centennial history of Hilbert's 16th problem, \emph{Bull. Am. Math. Soc} \textbf{39} (2002), 301-354 (new series).
		
		\bibitem{laf-1} G. Laferriere, G. J. Pappas and S. Yovine, Symbolic reachability computation for families of linear vector fields, \emph{J. Symb. Comput.} \textbf{32} (2001), 231-253.
		
		\bibitem{laf-2} G. Laferriere, G. J. Pappas and S. Yovine, A  new class of decidable hybrid systems, \emph{Hybrid Systems: Computation and Control  (LNCS vol 1569)} (1999) Berlin,  Springer, 137-151.	
		
		\bibitem{lli-mart-teix} J. Llibre, R. M. Martins and M. A. Teixeira, Periodic orbits, invariant tori and cylinders of Hamiltonian systems near integrable ones having a return map  equal to the  identity. \emph{J. Math. Phys} \textbf{51} (2010), 082704.
		
		\bibitem{lli-teix} J. Llibre and M. A. Teixeira, Limit cycles bifurcating from a $2$-dimensional isochronous cylinder, \emph{Appl. Math.  Lett.} \textbf{22} (2009), 321-341.
		
		\bibitem{lli-teix-torr} J. Llibre, M. A. Teixeira and J. Torregrosa, Limit cycles bifurcating form a $k$-dimensional isochronous center contained in $\mathbb{R}^n$ with $k \leq n$, \emph{Phys. Anal. Geom.} \textbf{10} (2007), 237-249.
		
		\bibitem{lli-zhang} J. Llibre and X. Zhang, Limit cycles of linear vector fields on  manifolds, \emph{Nonlinearity} \textbf{29} (2016), 3120-3131.
		
		\bibitem{malkin} I. G. Malkin, Some problems of the theory of nonlinear oscillations (Russian), \emph{Gosudarstv. Izdat. Tehn. - Theor. Lit.} (1956)
		
		\bibitem{rosseau} M. Rosseau, Vibrations non linéaires et théorie de la stabilité (French), \emph{Springer Tracts in Natural Philosophy}, vol 8 (1996) Berlin, Springer.
		
		
		\bibitem{san-ver-mur} J. A. Sanders, F. Verhulst and J. Murdock, Averaging Methods in Nonlinear Dynamical Systems \emph{Applied Mathematical Sciences)} vol 59 2nd edn (2007) New York, Springer.
		
		
		\bibitem{sotomayor} J. Sotomayor, Liçoes de equaçoes diferenciais ordinárias,  \emph{Euclid Project vol 11} (1979) Rio de Janeiro, Instituto de Matemática Pura e Aplicada.
		
		\bibitem{verhulst} F. Verhulst, Nonlinear Differential Equations and Dynamical Systems (1996) Berlin, Springer, Universitext.
		
	\end{thebibliography}
\end{document}